\documentclass[11pt]{article}
\usepackage{amsmath}
\usepackage{amssymb}
\usepackage{eucal}

\begin{document}

\baselineskip 16pt

\title{On ${\cal F}$-injectors of Fitting set of a finite group}

\author{Nanying Yang\thanks{Research of the  first  author is supported by
a NNSF grant of China (Grant \# 11301227)  and Natural Science Fund project in Jiangsu Province( grant \# BK20130119)}\\
{\small School of Science, Jiangnan University, Wuxi 214122, P. R. China}\\
 {\small E-mail: south0418@hotmail.com}\\ \\
W. Guo\thanks{Research of the second author is supported by
a NNSF grant of China (Grant \# 11371335)  and Wu Wen-Tsun Key Laboratory of Mathematics of Chinese Academy of Sciences.} \\
{\small Department of Mathematics, University
of Science and Technology of China,}\\ {\small Hefei 230026, P. R.
China}\\
 {\small E-mail: wbguo@ustc.edu.cn} \\ \\
{N.T. Vorob'ev}\\
{\small Department of Mathematics, Masherov Vitebsk State University,
Vitebsk 210038, Belarus}\\
 {\small E-mail: vorobyovnt@tut.by}}

\date{}
\maketitle

\begin{abstract}  Let $G$ be some generalized $\pi$-soluble groups and ${\cal F}$  be a Fitting set of $G$. In this paper, we prove the existence and conjugacy of ${\cal F}$-injectors of $G$, and give a description of the structure of the injectors.

\end{abstract}

\let\thefootnoteorig\thefootnote
\renewcommand{\thefootnote}{\empty}


\let\thefootnote\thefootnoteorig

\section{Introduction}

Throughout this paper, all groups are finite and $p$ is a prime. $G$ always denotes a group, $|G|$ is the order of $G$, $\sigma (G)$ is the set of all primes dividing $|G|$,  $\pi$ denotes a set of some primes. Let $\mathbb{P}$ be the set of all primes and $\pi'=\mathbb{P}\setminus \pi$. For any set ${\cal X}$ of subgroups of $G$, we let $\sigma ({\cal X})=\bigcup _{G\in {\cal X}}\sigma(G)$.

Recall that a class $\mathfrak{F}$ of groups is called a Fitting class if $\mathfrak{F}$ is closed under taking normal subgroups and products of normal $\frak{F}$-subgroups.
As usual, we denote by $\frak{E}, \frak{S}, \frak{N}$ the classes of all groups, all soluble groups, all nilpotent groups, respectively; $\frak{E}_{\pi}, \frak{S}_{\pi}, \frak{N}_{\pi}$  denote the classes of all $\pi$-groups, all soluble $\pi$-groups, all nilpotent $\pi$-groups, respectively; and $\frak{S}^{\pi}$ and $\frak{N}^{\pi}$ to denote the class of all $\pi$-soluble groups and the class of all $\pi$-nilpotent groups, respectively. It is well known that all the above classes are Fitting classes.

Following Anderson \cite{And} (see also \cite[VIII, (2.1)]{DH}), a nonempty set $\mathcal{F}$ of subgroups of $G$ is called a Fitting set of $G$
if the following three conditions are satisfied: i) If $T\unlhd S\in\mathcal{F}$, then  $T \in\mathcal{F}$; ii) If $ S,\ T \in\mathcal{F}$ and $S,\ T\unlhd ST$, then $ST\in\mathcal{F}$; iii) If $\ S\in\mathcal{F}$ and $x\in G$, then $S^x\in\mathcal{F}$.

For a set ${\cal X}$ of subgroups of $G$, the join $G_{\cal X}$ of all normal ${\cal X}$-subgroups of $G$ is called the ${\cal X}$-radical of $G$. From the definitions of Fitting class $\frak{F}$ and Fitting set ${\cal F}$ of $G$, for each nonempty Fitting class $\frak{F}$ (nonempty Fitting set ${\cal F}$ of $G$, resp.) of $G$, $G_{\frak{F}}$ ($G_{\cal F}$, resp.) is the unique maximal
normal $\frak{F}$-subgroup (${\cal F}$-subgroup resp.) of $G$.  In particular, $G_{\frak{N}}=F(G)$ is the Fitting subgroup of $G$.

It is well known that the theory of Fitting classes play an important role in the theory of groups. The importance of the theory of Fitting classes can firstly be seen in the following theorem proved by Fischer, Gasch$\ddot{u}$tz and Harley \cite{Fis1}: {\sl For every Fitting class $\frak{F}$, every  soluble group possesses exactly one conjugacy class of $\frak{F}$-injectors.} This
result is in fact a graceful generalization of the classical Sylow theorem and Hall theorem (The Hall Theorem shown that every soluble group possesses a Hall $\pi$-subgroup and any two Hall $\pi$-subgroups are conjugated in $G$ for any $\pi\subseteq \sigma (G)$).

Recall that, for any class $\frak{F}$ of groups, a subgroup $V$ of $G$ is said to be $\frak{F}$-maximal if $V\in \frak{F}$ and $U=V$ whenever $V\leq U\leq G$ and $U\in\frak{F}$.
A subgroup $V$ of a group $G$ is said to be an $\frak{F}$-injector of $G$ if $V\cap K$ is an
$\frak{F}$-maximal subgroup of $K$ for every subnormal subgroup $K$ of $G$. For a Fitting set  ${\cal F}$ of $G$, the ${\cal F}$-injector of $G$ is similarly defined (see \cite[VIII, (2.5)]{DH}).
Clearly, every $\frak{F}$-injector (resp. ${\cal F}$-injector) of $G$ is an $\frak{F}$-maximal (resp. ${\cal F}$-maximal) subgroup of $G$.
Note that if $\frak{F}=\frak{N}_p$ is the Fitting class of all $p$-groups, then the $\frak{F}$-injectors of a group $G$ are Sylow $p$-subgroups of $G$; if $\frak{F}=\frak{E}_{\pi}$  and $G$ has a Hall $\pi$-soluble group, then the $\frak{F}$-injectors of $G$ are Hall $\pi$-subgroups of $G$ (see \cite[p.68, Ex.1]{Guo1} or \cite[p.328]{BalE}).

As a development of the theorem of Fischer, Gasch$\ddot{u}$tz and Harley \cite{Fis1},  Shemetkov \cite{Shem1} (Anderson \cite{And}) proved that if $G$ is a $\pi$-soluble group (resp. soluble group) and  ${\cal F}$ is a Fitting set of $G$, then $G$ possesses exactly one conjugacy class of ${\cal F}$-injectors, where $\pi =\sigma({\cal F})$.

If $\frak{F}$ is a Fitting class, then the set $\{H\leq G :  H\in \frak{F}\}$ is a Fitting set, which is denoted by $Tr_{\frak{F}}(G)$ and called the trace of $\frak{F}$ in $G$ (see \cite[p.537]{DH}). Note that for a Fitting class $\frak{F}$, the $\frak{F}$-injectors and $Tr_{\frak{F}}(G)$-injectors of $G$ coincide (see below Lemma 2.1).  But by \cite[VIII, Examples (2.2)(c)]{DH}, not every  Fitting set of $G$ is the trace of a Fitting class. Hence, in view of the theorem of Shemetkov in \cite{Shem1} and the theorem of  Fischer, Gasch$\ddot{u}$tz and Harley \cite{Fis1}, the following  question  naturally arise:

{\bf Problem 1.1.} For an arbitrary Fitting set ${\cal F}$ of a non-soluble group $G$ (in particular, some generalized $\pi$-soluble groups $G$, for example, the groups in Theorem A below), whether $G$ possesses an ${\cal F}$-injector and any two ${\cal F}$-injectors are conjugate ?

It is clear that for a Fitting class $\frak{F}$, every $\frak{F}$-injector of $G$  contains the $\frak{F}$-radical $G_{\frak{F}}$ of $G$. In \cite{Fischer} Fischer gave the characterization the $\frak{N}$-injectors by means of the nilpotent radical $G_{\frak{N}}=F(G)$. He proved that the set of $\frak{N}$-injectors is exactly the set of all maximal nilpotent subgroups of $G$ containing $F(G).$

The product $\frak{FH}$ of two Fitting classes $\frak{F}$ and
$\frak{H}$ is the class $(G \mid G/G_{\frak{F}}\in \frak{H}).$ It is
well known that the product of any two Fitting classes is also a
Fitting class and the multiplication of Fitting classes satisfies
associative law (see \cite[IX, (1.12)(a)(c)]{DH}).

Hartley \cite{Hart} proved that, for the Fitting class of type $\frak{XN}$ (where $\frak{X}$ is a nonempty Fitting class), a subgroup $V$ of a soluble group $G$ is an $\frak{XN}$-injecotr of $G$ if and only if $V/G_{\frak{X}}$ is a nilpotent subgroup of $G$. As a further improvement, the authors in \cite{Guo137} proved that for a Hartley class $\frak{H}$, a subgroup $V$ of soluble group $G$ is $\frak{H}$-injector of $G$ if and only if $V/G_h$ is a $D$-injector of $G/G_h$, where $D=\bigcap _{i\in I}\frak{S}_{\pi '}\frak{S}_{\pi}$ and $G_h=\prod _{i\in I}G_{h(\pi _i)}$. Moreover, in \cite{Guo137}, it was proved that the set of $\frak{H}$-injectors coincide with the set of of all $\frak{H}$-maximal subgroups containing the $\frak{H}$-radical of $G$.

Following \cite[\S 7.2]{BalE}, a Fitting set ${\cal F}$ is said to be  injective if every group $G$ possesses at least one ${\cal F}$-injector.

In connection with above, the following more general question naturally arise.

{\bf Problem 1.2.} Let ${\cal F}$ be a injective Fitting set of $G$ in some universe. What is the structure of the ${\cal F}$-injector?

In order to resolve the Problems 1.1 and 1.2, we need to develop and extend the local method of Hartley \cite{Hart} (which is for Fitting classes and in the universe of soluble groups) to for Fitting sets and in a more general universe (not necessary in the universe of soluble groups).

Firstly, for a Fitting set ${\cal F}$ of $G$ and a Fitting class $\frak{X}$, we call the set $\{H\leq G | H/H_{\cal F}\in \frak{X}\}$ of subgroups of $G$ the product of ${\cal F}$ and $\frak{X}$ and denote it by ${\cal F}\circ \frak{X}.$

Following \cite{Vorb1}, a function $f$ : $\mathbb{P} \longrightarrow \{\rm{Fitting \ sets \ of \ G}\}$ is called a Hartley function (or in brevity, H-function).

{\bf Definition 1.3.} Let $\emptyset\neq \pi \subseteq \mathbb{P}$ and
$SLR(f)=\bigcap _{p\in \pi}f(p)\circ \frak{E}_{p'}.$  A Fitting set ${\cal F}$ of $G$ called: {\sl $\pi$-semilocal} if ${\cal F}=SLR(f)$ for some $H$-function $f$; {\sl $\pi$-local}
if ${\cal F}=SLR(f)$ for some $H$-function $f$ such that $f(p)\circ\frak{N}_p=f(p)$ for all $p\in \pi$. In this case when ${\cal F}=SLR(f)$, ${\cal F}$ is said to be defined by the H-function $f$ or $f$ is an H-function of ${\cal F}$.

In particular, if $\pi=\mathbb{P}$, then a  $\pi$-semilocal Fitting set (a $\pi$-local Fitting set, respectively) is said to be a semilocal Fitting set (a local Fitting set, respectively).

{\bf Definition 1.4.} Let ${\cal F}=SLR(f)$ be a $\pi$-semilocal Fitting set of $G$.  Then $f$ is said to be

1)  integrated if $f(p)\subseteq {\cal F}$ for all $p\in \pi$;

2) full if $f(p)=f(p)\circ\frak{N}_p$ for all  $p\in \pi$;

3) invariable if $f(p)=f(q)$ for all $p, q\in \pi$.

For a $\pi$-semilocal Fitting set ${\cal F}=SLR(f)$ of $G$, the subgroup $G_f=\Pi _{p\in \pi}G_{f(p)}$ is said the $f$-radical of $G$.

The following two theorems resolved the Problems 1.1 and 1.2 for some generalized $\pi$-soluble groups  and some $\pi$-semilocal Fitting sets of $G$.

{\bf Theorem A. } Let ${\cal F}$ be a Fitting set of $G$. Then $G$ possesses an ${\cal F}$-injector and any two ${\cal F}$-injectors are conjugate if one of the following conditions:

(1) $G\in {\cal F}\circ\frak{S}^{\pi}$, where $\pi =\sigma ({\cal F});$

(2) ${\cal F}$ is $\pi$-semilocal and $ G\in \frak{S}^{\pi};$

(3) ${\cal F}$ is $\pi$-semilocal and $G\in {\cal F}\circ\frak{S}^{\pi}$, where $\sigma(G_{\cal F})\subseteq \pi $.

Moreover, the index of every ${\cal F}$-injector in $G$ is a $\pi$-number in the case (2).

{\bf Theorem B.} Let ${\cal F}$ be a $\pi$-semilocal Fitting set of $G$ defined by a full and invariable H-function $f$ and $G_f$ the $f$-radical of $G$.
If $G/G_f$ is $\pi$-soluble, then the following statements hold:

(1) A subgroup $V$ of $G$ is an ${\cal F}$-injector of $G$ if and only if $V/G_f$ is a Hall $\pi '$-subgroup of $G/G_f;$

(2) $G$ possesses an ${\cal F}$-injector and any two ${\cal F}$-injectors are conjugate in $G$;

(3) A subgroup $V$ of $G$  is an ${\cal F}$-injector of $G$ if and only if $V$ is an ${\cal F}$-maximal subgroup of $G$ and $G_{{\cal F}}\leq V;$

(4) If $G\in \frak{S}^{\pi}$, then every ${\cal F}$-injector of $G$ is of the type $G_{\pi '}G_f$, where $G_{\pi '}$ is some Hall $\pi '$-subgroup of $G$.

Theorem A give the new theory of ${\cal F}$-injectors for the generalized $\pi$-soluble group (in particular, $\pi$-soluble group). Theorem B describe the structures of the injectors. From Theorems A and B, a series of famous results can be directly follow. For example, Fischer, Gasch$\ddot{u}$tz and Harley Theorem \cite{Fis1} (see also \cite[VIII, Theorem (2.9) and IX, Theorem (1.4)]{DH}),
Shemetkov \cite[Theorem 2.2]{Shem1}, Ballester-Bolinches and Ezquerro \cite[Theorem 2.4.27]{BalE}, Sementovskii \cite[Theorem]{Semen}, Bolado-Caballero, Martiner-Verduch \cite[Theorem 2.2]{Bolado} and Guo \cite{Guo2} (see also \cite[Theorem 2.5.3]{Guo1}.

{\bf Remark 1.5.} (a) The statement of Theorem A is not true in general for any Fitting set and a non-soluble group, for example, let $S=A_7$ the alternating group of degree $7$, $T=PSL(2, 11)$ the projective special linear group, a Fitting class $\frak{F}=D_0(S, T, 1)$ and a group $G$ such as in \cite[Theorem 7.1.3]{BalE}. Then for the Fitting set  ${\cal F}=Tr_{\frak{F}}(G)$, $G$ has no ${\cal F}$-injector (see \cite[Theorem 7.1.3]{BalE}).

(b) It is well known that there exist examples of the sets $\pi$ of primes and non-abelian groups $G$ such that $G$ possesses a Hall $\pi$-subgroup but the Hall $\pi$-subgroups are not conjugate (see Revin and Vdovin \cite[Theorem 1.1]{RevinV}). However, a Hall $\pi$-subgroups of a $E_{\pi}$-group $G$ (Note that a group $G$ is said to be an $E_{\pi}$-group if $G$ has a Hall $\pi$-subgroup) is a $\frak{E}_{\pi}$-injector of $G$ (see \cite[p. 328]{BalE}). Therefore,  $\frak{E}_{\pi}$-injectors of a $E_{\pi}$-group $G$ are not conjugate in general.

The following example shows that the statement of Theorem B(3) is not true if  ${\cal F}$ is only a Fitting set even of a soluble $G$.

{\bf Example 1.6.} Let a Fitting class $\frak{F}=\frak{N}\frak{N}_3$, $S=S_3$ the symmetric group of degree $3$ and $M$ a faithful irreducible $S$-module over $F_5$,  and let $T=M\rtimes S$, $N$ be a faithful irreducible $T$-module over $F_2$, $G=N\rtimes T$  and ${\cal F}=Tr_{\frak{F}}(G)$. Then $F(G)=N$ is a 2-group, and clearly  ${\cal F}$ is not $\pi$-semilocal for any $\pi\subsetneq \mathbb{P}$. Similar as \cite[IX, Example 4.4]{DH}, we know that the set of ${\cal F}$-injectors coincide with the set of the subgroups of type $F(G)P$, where $P$ is some Sylow 3-subgroup of $G$. But every Sylow 2-subgroup of $G$ is an ${\cal F}$-subgroup of $G$ containing the ${\cal F}$-radical $G_{{\cal F}}=F(G)$ and it is not contained in any ${\cal F}$-injector of $G$.

All unexplained notion and terminology are standard. The reader is referred to \cite{DH, BalE, Guob}.

\section{Preliminaries}

 \ \ \ \ {\bf Lemma 2.1} (see \cite[VIII, (2.4)(d) and p. 601]{DH}).  1) Let ${\cal F}$ be a Fitting set of $G$. If $N$ is a subnormal subgroup of $G$, then $N_{{\cal F}}=N\cap G_{{\cal F}}.$

 2) If $\frak{F}$ is a Fitting class, then the $\frak{F}$-injectors and $Tr_{\frak{F}}(G)$-injectors of $G$ coincide.

In view of the definition of ${\cal F}$-injector and \cite[IX, (1.3), VIII, (2.6)(2.7)]{DH}, we have the following

{\bf Lemma 2.2.}  Let ${\cal F}$ be a Fitting set of a groups $G$.

(1) If $V$ is an ${\cal F}$-injector  of $G$, then $G_{\cal F}\leq V;$

(2) If $V$ is an ${\cal F}$-injector of $G$, then  $V$ is an ${\cal F}$-maximal subgroup of $G$.

(3) If $V$ is an ${\cal F}$-maximal subgroup of $G$ and $V\cap M$ is an ${\cal F}$-injector  of $M$ for any maximal normal subgroup $M$ of $G$, then $V$ is an ${\cal F}$-injector  of $G$;

(4) If $N\unlhd G$ and $V$ is an ${\cal F}$-injector  of $N$, then $V^g$ is an ${\cal F}$-injector  of $N$ for all $g\in G;$

(5) If $K$ is a subnormal subgroup of $G$ and  $V$ is an ${\cal F}$-injector of $G$, then $V\cap K$ is an  ${\cal F}$-injector of $K$.

{\bf Lemma 2.3} \cite{Shem1}.   Let ${\cal F}$ be a Fitting set of $G$ and $\pi=\sigma (\cal F)$. If $G$ is $\pi$-soluble, then $G$ has an ${\cal F}$-injector and any two ${\cal F}$-injector of $G$ are conjugate in $G$.

{\bf Lemma 2.4} \cite[Lemma 3.9]{Vorb2}.  Let $G$ be a $\pi$-soluble group and ${\cal F}$ a Fitting set of $G$ such that ${\cal F}\circ \frak{E}_{\pi'}={\cal F}$. Assume that $N\unlhd G$ and $G/N$ is either a nilpotent $\pi$-group or $\pi'$-soluble. Let $W$ is an ${\cal F}$-maximal subgroup of $N$. If $V_1$ and $V_2$ are ${\cal F}$-maximal subgroup of $G$ such that $W\leq V_1\cap V_2$, then $V_1$ and $V_2$ are conjugate in $G$.

{\bf Lemma 2.5} \cite[Theorems 2.4.27]{BalE}. Let ${\cal F}$ be a Fitting set of $G$ and $G/G_{\cal F}$ is soluble. Then $G$ possesses exactly one conjugacy class of ${\cal F}$-injectors.

\section{Properties of the product of a Fitting set and a Fitting class}

In this section, we discuss the properties of the product of a Fitting set and a Fitting class, which are also useful in the proof of Theorems A and B.

{\bf Proposition 3.1.} Support that ${\cal F}$ is a Fitting set of $G$ and $\frak{H}$ is a nonempty Fitting class. Then ${\cal F}\circ \frak{H}$ is a Fitting set of $G$.

{\bf Proof.} Let ${\cal M}={\cal F}\circ \frak{H}$. Assume that $H\leq G$, $K\unlhd H\in {\cal M}.$ Then $H/H_{{\cal F}}\in \frak{H}$. By Lemma 2.1, $K_{{\cal F}}=K\cap H_{{\cal F}}.$  Hence $K/K_{{\cal F}}\simeq KH_{{\cal F}}/H_{{\cal F}}\unlhd H/H_{{\cal F}}\in \frak{H}$. Consequently, $K/K_{{\cal F}}\in \frak{H}$, and so $K\in {\cal M}.$

Now assume that $H, L$ are subgroup of $G$, $H, L\in {\cal M}$ and $H, L\unlhd HL$. Then $H/H_{\cal F}\in \frak{H}$. By Lemma 2.1, $H_{\cal F}=(HL)_{\cal F}\cap H$. Similarly,  $L_{\cal F}=(HL)_{\cal F}\cap L$. By the isomorphism
$H/H_{\cal F}\simeq H(HL)_{\cal F}/(HL)_{\cal F}$  and $L/L_{\cal F}\simeq L(HL)_{\cal F}/(HL)_{\cal F}$, we have that
$HL/(HL)_{\cal F}\in \frak{H}.$  Hence $HL\in {\cal M}.$

 Final assume that $H\leq G$ and $H\in {\cal M}.$ Then $\bar{H}=H/H_{{\cal F}}\in \frak{H}$ and $\bar{H}=(\bar{H})^x =H^x/(H_{{\cal F}})^x.$ In order to prove that $H^x\in {\cal M}$ for all $x\in G$, it only need to prove that $(H^x)_{{\cal F}}=(H_{{\cal F}})^x.$
In fact, since $H_{{\cal F}}\in {\cal F}$ and ${\cal F}$ is a Fitting set of $G$, $(H_{{\cal F}})^x \in {\cal F}$. Moreover, $(H_{{\cal F}})^x\unlhd H^x$,  so  $(H_{{\cal F}})^x\leq (H^x)_{{\cal F}}.$ On the other hand, $(H^x)_{{\cal F}}\unlhd H^x$, so  $((H^x)_{{\cal F}})^{x^{-1}}\unlhd H$. But as ${\cal F}$ is a Fitting set of $G$ and $(H^x)_{{\cal F}}\in {\cal F},$ we have that $((H^x)_{{\cal F}})^{x^{-1}}\in {\cal F}$, and so  $((H^x)_{{\cal F}})^{x^{-1}}\leq H_{{\cal F}}.$  It follows that $(H^x)_{{\cal F}}=(H_{{\cal F}})^x.$

The above shows that ${\cal F}\circ \frak{H}$ is a Fitting set of $G$. The proposition is proved.

{\bf Proposition 3.2.} Let $H\leq G$ and ${\cal F}$ be a Fitting set of $G$, and $\frak{H}$ a nonempty Fitting class. Then:

(1) ${\cal F}\subseteq {\cal F}\circ \frak{H}.$

(2) $H_{{\cal F}\circ \frak{H}}/H_{{\cal F}}=(H/H_{{\cal F}})_{\frak{H}}.$

{\bf Proof.} (1) Assume that $H\in {\cal F}$. Then $H=H_{{\cal F}},$ so $H/H_{{\cal F}}=1\in \frak{H}$. It follows that $H\in {\cal F}\circ \frak{H}.$ Hence (1) holds.

(2) Since ${\cal F}\subseteq {\cal F}\circ \frak{H}$ by (1), $H_{{\cal F}}\leq H_{{\cal F}\circ \frak{H}}\leq H.$ Then by Lemma 2.1(1), $H_{{\cal F}}=H_{{\cal F}\circ \frak{H}}\cap H_{\cal F}= (H_{{\cal F}\circ \frak{H}})_{{\cal F}}. $
As $H_{{\cal F}\circ \frak{H}}\in {\cal F}\circ \frak{H},$ $H_{{\cal F}\circ \frak{H}}/(H_{{\cal F}\circ \frak{H}})_{{\cal F}}\in \frak{H}.$ Hence $H_{{\cal F}\circ \frak{H}}/H_{{\cal F}}\leq (H/H_{{\cal F}})_{\frak{H}}.$

Now let $L/H_{{\cal F}}=(H/H_{{\cal F}})_{\frak{H}}.$
Then $L/H_{{\cal F}}\unlhd H/H_{{\cal F}},$  so $L\unlhd H$. By Lemma 2.1, $L_{{\cal F}}=L\cap H_{{\cal F}}=H_{{\cal F}}.$ Hence $L/H_{{\cal F}}=L/L_{{\cal F}}\in \frak{H}$. Consequently, $L\in {\cal F}\circ \frak{H}.$ But as $L\unlhd H$, we obtain that $L\leq H_{{\cal F}\circ \frak{H}}.$ Thus $L/H_{{\cal F}}\leq H_{{\cal F}\circ \frak{H}}/H_{{\cal F}}.$

{\bf Proposition 3.3.} Let ${\cal F}$ be a Fitting set of $G$, $H\leq G$ and $\frak{M, H}$ be nonempty Fitting classes. Then
$$({\cal F}\circ \frak{M})\circ \frak{H}={\cal F}\circ \frak{MH}.$$

{\bf Proof.} Let $H\in ({\cal F}\circ \frak{M})\circ \frak{H}$. Then $H/H_{{\cal F}\circ \frak{M}}\in \frak{H}.$ Hence by Proposition 3.2,
$$H/H_{{\cal F}\circ \frak{M}}\simeq (H/H_{{\cal F}})/(H_{{\cal F}\circ \frak{M}}/H_{{\cal F}})=(H/H_{{\cal F}})/(H/H_{{\cal F}})_{\frak{M}}\in \frak{H}.$$
Hence $H/H_{{\cal F}}\in \frak{MH}$, and so $H\in {\cal F}\circ (\frak{MH}).$ This implies that $({\cal F}\circ \frak{M})\circ \frak{H}\subseteq {\cal F}\circ \frak{MH}.$

Since the above derivation process is reversible, we have also that ${\cal F}\circ \frak{MH}\subseteq ({\cal F}\circ \frak{M})\circ \frak{H}.$  Thus the proposition holds.

Recall that a class $\frak{X}$ of groups is called a formation if it is closed under taking homomorphic images and subdirect products. Clearly, for a nonempty formation, every group $G$ has the least normal subgroup $G^{\frak{X}}$ such that $G/G^{\frak{X}}\in \frak{X}.$  A Fitting formation is a class of groups which is both a formation and a Fitting class.  A Fitting class which is closed under homomorphic image is called a radical homomorph (see, for example, \cite[p. 2 and p. 324]{Guob}).

{\bf Proposition 3.4.} Let ${\cal F, H}$ be two Fitting sets of $G$. Then

(1) If $\frak{M}$ is a nonempty radical homomorph and  ${\cal F}\subseteq {\cal H},$ then ${\cal F}\circ \frak{M}\subseteq {\cal H}\circ \frak{M}.$

(2) If $\frak{M}$ is a nonempty Fitting formation, then $({\cal F}\cap {\cal H})\circ\frak{M}=({\cal F}\circ\frak{M})\cap ({\cal H}\circ\frak{M}).$

(3) If $\frak{M}$  and $\frak{H}$ are nonempty Fitting classes, then ${\cal F}\circ (\frak{M}\cap \frak{H})=({\cal F}\circ\frak{M})\cap ({\cal F}\circ\frak{H}).$

{\bf Proof.} The proof is similar to the proof of \cite[Lemma 4]{Vorb3}.

\section{Semilocal Fitting set and $f$-radical}

In this section, we give some results about semilocal Fitting sets and $f$-radical, which are also main steps in the proof of Theorems A and B.

Firstly, we give some following examples of semilocal Fitting sets.

{\bf Example 4.1.} Support that ${\cal F}$ is a Fitting set of $G$ and  $\pi\subseteq \mathbb{P}$.

(a) If an H-function $f$ satisfy $f(p)=\{1\}$ for all $p\in \pi$, and ${\cal F}=SLR(f)$, then by definition of $SLR(f)$, ${\cal F}$ is the set of all $\pi '$-subgroups of $G$. In particular, if $\pi =\mathbb{P}$, then ${\cal F}=\{1\}.$

(b) If $f(p)=\{1\}\circ \frak{N}_p$ for all $p\in \pi$, then by Proposition 3.3, $f(p)=(\{1\}\circ \frak{N}_p)\circ \frak{E}_{p'}=\{1\}\circ (\frak{N}_p\frak{E}_{p'})$  for all $p\in \pi$. If ${\cal F}=SLR(f)$, then ${\cal F}$ is the set of all the subgroups of $G$ which have normal nilpotent Hall $\pi$-subgroup (in fact, if $H\in SLR(f)$, then $H\simeq H/H_{\{1\}}\in \bigcap _{p\in \pi}\frak{N}_p\frak{E}_{p'}=\frak{N}_{\pi}\frak{E}_{\pi '}$).  In particular, if $\pi =\mathbb{P},$ then  ${\cal F}$ is the set of all nilpotent subgroups of $G$.

(c) Assume that $G$ is $\pi$-soluble and ${\cal M}=\{H\leq G | H_{\pi}\in{\cal F}\}.$ Then $G$ is a $D_{\pi}$-group, that is, $G$ has a Hall $\pi$-subgroup and every $\pi$-subgroup of $G$ is conjugately contained some Hall $\pi$-subgroup of $G$ (see \cite{Cun4} or \cite[Theorem 1.7.6]{Guo1}. Since ${\cal F}$ is a Fitting set, it is easy to see that ${\cal M}$ is a Fitting set of $G$.
Now we prove that ${\cal M}$ is $\pi$-semilocal. By Proposition 3.2(1), ${\cal M}\subseteq {\cal M}\circ \frak{E}_{\pi '}.$ Let $H\leq G$ and $H\in {\cal M}\circ \frak{E}_{\pi '}.$
Then $H_{{\cal M}}\in \frak{E}_{\pi'}$, and so $H_{\pi}\leq H_{{\cal M}}.$ Hence  $H_{\pi}\leq (H_{{\cal M}})_{\pi}=H_{\pi} \cap H_{{\cal M}}\leq H_{\pi}.$ This implies that $H_{\pi}=(H_{{\cal M}})_{\pi}$. Since $H_{{\cal M}}\in {\cal M}$, $H_{\pi}\in {\cal F}.$ Hence $H\in{\cal M}$. This shows that ${\cal M}\circ \frak{E}_{\pi '}\subseteq {\cal M}.$ Therefore ${\cal M}\circ \frak{E}_{\pi '}= {\cal M}.$
Now by Proposition 3.4(3),  ${\cal M}={\cal M}\circ \frak{E}_{\pi '}={\cal M}\circ (\bigcap_{p\in \pi}\frak{E}_{p'})=\bigcap_{p\in \pi}{\cal M}\circ \frak{E}_{p'}=SLR(f)$  for the H-function $F$ such that $f(p)={\cal M}$ for all $p\in \pi.$  Therefore, ${\cal M}$ is semilocal.

The following result give the characterization of $\pi$-semilocal Fitting set.

{\bf Lemma 4.2.} Let ${\cal F}$ be a Fitting set of $G$. Then ${\cal F}$ is $\pi$-semilocal if and only if ${\cal F}\circ \frak{E}_{\pi '}={\cal F}.$

{\bf Proof.} Assume that ${\cal F}$ is $\pi$-semilocal. If $\pi=\mathbb{P},$ then the lemma is clear. We may, therefore, assume that $\pi\neq \mathbb{P}.$  Then ${\cal F}\circ \frak{E}_{\pi '}=(\bigcap _{p\in \pi}f(p)\circ \frak{E}_{p'})\circ \frak{E}_{\pi '}.$  By Propositions 3.3 and 3.4(2),
$${\cal F}\circ \frak{E}_{\pi '}=\bigcap _{p\in \pi}f(p)\circ (\frak{E}_{p'}\frak{E}_{\pi '})=\bigcap _{p\in \pi}f(p)\circ \frak{E}_{p'}={\cal F}.$$

Conversely, if ${\cal F}\circ \frak{E}_{\pi '}={\cal F}$, then ${\cal F}={\cal F}\circ (\bigcap _{p\in \pi}\frak{E}_{p '}).$ By Proposition 3.4(3), ${\cal F}=\bigcap _{p\in \pi}{\cal F}\circ \frak{E}_{p'}$. Hence ${\cal F}=SLR(f)$ for the H-function $f$ such that $f(p)={\cal F}$ for all $p\in \pi.$  This completes the proof.

For an H-function $f$ of a semilocal Fitting set of $G$, we call the subgroup $G_f=\Pi_{p\in \pi}G_{f(p)}$ of $G$ the {\sl $f$-radical} of $G$.

{\bf Lemma 4.3.} Support that ${\cal F}$ is a $\pi$-semilocal Fitting set of $G$ defined by a integrated H-function $f$. If $H\leq G$ and $H/G_f\in \frak{E}_{\pi '}$, then $H\in{\cal F}.$

{\bf Proof.} Let $p\in \pi$. Clearly, $G_{f(p)}=(G_f)_{f(p)}=G_f\cap H_{f(p)}$  by Lemma 2.1(1). Hence $G_{f(p)}\leq H_{f(p)}.$  Since $f$ is integrated, $G_f\in {\cal F}.$ It follows that $G_f/G_{f(p)}\in \frak{E}_{p'}.$ Then by the isomorphism
$$H_{f(p)}G_f/H_{f(p)}\simeq G_f/(G_f\cap H_{f(p)})\simeq (G_f/G_{f(p)})/((G_f\cap H_{f(p)})/G_{f(p)}),$$
we have that $H_{f(p)}G_f/H_{f(p)}\in \frak{E}_{p'}.$ Since $H/G_f\in \frak{E}_{\pi '}=\bigcap _{p\in \pi}\frak{E}_{p'}$ by the hypothesis, by the isomorphism
$H/H_{f(p)}G_f\simeq (H/G_f)/(H_{f(p)}G_f/G_f),$ we have also that
$H/H_{f(p)}G_f\in \frak{E}_{\pi '}\subseteq \frak{E}_{p'}.$ This implies that $H/H_{f(p)}\in \frak{E}_{p'}$ for all $p\in \pi$. Thus $H\in\bigcap _{p\in \pi}f(p)\circ\frak{E}_{p'}={\cal F}.$

{\bf Lemma 4.4.} Support that ${\cal F}$ is a $\pi$-semilocal Fitting set of $G$ defined by a full H-function $f$. If $G/G_{\cal F}$ is $\pi$-soluble and ${\cal X}$ is a Fitting set of $G$ such that ${\cal X}\subseteq \bigcap _{p\in \pi}f(p)$, then $C_G(G_{{\cal F}}/G_{{\cal X}})\leq G_{{\cal F}}.$

{\bf Proof.} Suppose that this theorem is false and let $G$ be a counterexample of minimal order. Then clearly, $|G|\neq 1.$  Let $C=C_G(G_{{\cal F}}/G_{{\cal X}}).$ Then by $C\nleq G_{{\cal F}}$, we have that $C/C\cap G_{{\cal F}}\neq 1.$ Hence there exists a normal subgroup $N$ of $G$ such that $N/N\cap G_{{\cal F}}=N/C\cap G_{{\cal F}}$ is a nontrivial chief factor of $G$ and $G_{\frak{X}}\leq V.$ Since $N/N\cap G_{{\cal F}}\simeq NG_{{\cal F}}/G_{{\cal F}}$ and $G\in {\cal F}\circ \frak{S}^{\pi},$ we have that $N/N\cap G_{{\cal F}}$ is a $\pi$-soluble group. Hence $N/N\cap G_{{\cal F}}$ is either a $\pi'$-group or an elementary abelian $p$-group for some $p\in \pi.$

Assume that $N/N\cap G_{{\cal F}}\in \frak{E}_{\pi '}$. Then by Lemma 2.1(1), $N_{\cal F}=N\cap G_{\cal F},$ so $N/N_{{\cal F}}\in \frak{E}_{\pi '}.$ It follows from Lemma 4.2 that $N\in {\cal F}\circ \frak{E}_{\pi '}={\cal F}$. Hence $N\leq G_{{\cal F}},$ a contradiction.

Now assume that $N/N\cap G_{{\cal F}}\in \frak{A}_p$, where $\frak{A}_p$ is the class of all abelian $p$-groups. Then the derived subgroup $(N/N\cap G_{{\cal F}})'=1$. Hence $N'(N\cap G_{{\cal F}})=N\cap G_{{\cal F}}$,  and so $N'\leq N\cap G_{{\cal F}}.$ As $N\leq C$, we have that $N\leq C_G(N\cap G_{{\cal F}}/G_{{\cal X}}).$
Hence $[N', N]\leq [N\cap G_{{\cal F}}, N]\leq G_{{\cal X}}.$  This implies that $[(N/G_{{\cal X}})', N/G_{{\cal X}}]=1$, so $N/G_{{\cal X}}$ is nilpotent with nilpotent class $\leq 2$.
Let $\bar{P}=L/G_{{\cal X}}$ be the normal Sylow $p$-subgroup of $N/G_{{\cal X}}$. Then $L\unlhd G$. By \cite[A, (10.9)]{DH}, $\bar{P}$ covers the chief factor
$(N/G_{{\cal X}})/(N\cap G_{{\cal X}}/G_{{\cal X}})$ of $G/G_{{\cal X}}$, that is, $N/G_{{\cal X}}\leq \bar{P}(N\cap G_{{\cal X}}/G_{{\cal X}})$. Hence $L(N\cap G_{{\cal F}})\geq N$.
It follows that $LG_{{\cal F}}=NG_{{\cal F}}.$
Since $L_{{\cal X}}=L\cap G_{{\cal X}}=G_{\cal X}$ by Lemma 2.1, $L\in {\cal X}\circ \frak{N}_p\subseteq (f(p)\circ\frak{N}_p)\circ \frak{E}_{p'}=f(p)\circ(\frak{N}_p\frak{E}_{p'})$ by Proposition 3.3. Let $p, q\in \pi$ and $p\neq q$. By Propositions 3.3 and 3.4, $L\in {\cal X}\circ \frak{N}_p\subseteq f(q)\circ\frak{E}_{q'} \subseteq (f(q)\circ \frak{N}_{q})\circ \frak{E}_{q'}=f(q)\circ(\frak{N}_q\frak{E}_{q'})$. Note that $f$ is full, so $f(p)\circ \frak{N}_p=f(p)$. Hence  $L\in \bigcap _{p\in \pi}f(p)\circ(\frak{N}_p\frak{E}_{p'})=\bigcap_{p\in \pi}(f(p)\circ\frak{N}_p)\circ \frak{E}_{p'}=\bigcap _{p\in\pi}f(p)\circ\frak{E}_{p'}={\cal F}.$   Consequently, $L\leq G_{{\cal F}}$, and so
$NG_{{\cal F}}=LG_{\cal F}\leq G_{{\cal F}},$ which contradicts that $NG_{{\cal F}}/G_{{\cal F}}$ is a chief factor of $G$.
This contradiction completes the proof.

A group $G$ is said to be $\frak{F}$-constrained if $C_G(G_{\frak{F}})\subseteq G_{\frak{F}}.$ It is well known that if $G$ is soluble, then $G$ is $\frak{N}$-constrained, that is, $C_G(F(G))\subseteq F(G).$

Recall that $G$ is said to be: $\pi$-closed if $G$ has a normal Hall $\pi$-subgroup; $\pi$-special if $G$ has a normal nilpotent Hall $\pi$-subgroup (see \cite{Cun3}). Let $\frak{F}$ is the class of all $\pi$-closed groups and $\frak{H}$ be the class of all $\pi$-special groups, that is, $\frak{F}=\frak{E}_{\pi}\frak{E}_{\pi '}$ and $\frak{H}=\frak{N}_{\pi}\frak{E}_{\pi '}.$ Then $G_{\frak{F}}$ is called the $\pi$-closed radical of $G$, and  $G_{\frak{H}}$ is called the $\pi$-special radical of $G$.
Let ${\cal F}=Tr_{\frak{F}}(G)$ and ${\cal H}=Tr_{\frak{H}}(G)$. Then ${\cal F}$ is the set of all $\pi $-closed subgroups of $G$ and  ${\cal H}$ is the set of all $\pi$-special subgroups of $G$. Obviously, $G_{\cal F}=G_{\frak{F}}$ and $G_{\cal H}=G_{\frak{H}}$ (see \cite[p.563]{DH}).  Hence by Lemma 4.4, we obtain that $G$ is $\frak{F}$-constrained and $\frak{H}$-constrained (it only need let $\frak{X}=\{1\}$ in Lemma 4.4).

{\bf Corollary 4.5.} If $\frak{F}=\frak{E}_{\pi}\frak{E}_{\pi '}$, then every $\pi$-soluble group is $\frak{F}$-constrained.

{\bf Proof.}  Let $G$ be a $\pi$-soluble group and ${\cal F} =Tr_{\frak{F}}(G)$.

We first prove that ${\cal F}=SLR(f)$ for the full H-function $f$ such that $f(p)=Tr_{\frak{E}_{\pi}}(G)$ for all $p\in \pi$.

In fact, if $X\leq G$ and $X\in f(p)\circ \frak{N}_p$ for all $p\in \pi$, then $X/X_{f(p)}\in \frak{E}_{p}$ and so $X\in \frak{E}_{\pi}$. Hence $X\in Tr_{\frak{E}_{\pi}}(G)=f(p).$ This shows that $f(p)\circ \frak{N}_p\subseteq f(p).$ But by Proposition 3.2, $f(p)\subseteq f(p)\circ \frak{N}_p.$ Thus, $f(p)=f(p)\circ \frak{N}_p$, this is, $f$ is full.

Let $L\leq G$ and $L\in {\cal F}$. Then $L\in \frak{E}_{\pi}\frak{E}_{\pi '}$, so $L/L_{\frak{E}_{\pi}}\in\frak{E}_{\pi '}.$
But, clearly,  $L_{\frak{E}_{\pi}}\leq L_{f(p)}$ for all
$p\in \pi$. By the isomorphism $(L/L_{\frak{E}_{\pi}})/(L_{f(p)}/L_{\frak{E}_{\pi}})\simeq L/L_{f(p)},$ we have that $L/L_{f(p)}\in \frak{E}_{\pi '}\subseteq \frak{E}_{p'}$ for all $p\in \pi$. This implies that $L\in \bigcap _{p\in \pi}f(p)\circ \frak{E}_{p'}=SLR(f)$. Hence ${\cal F}\subseteq SLR(f).$

Now assume that $R\in SLR(f).$ Then $R\in f(p)\circ \frak{E}_{p'}$ and so $R/R_{f(p)}\in \frak{E}_{p'}$ for all $p\in \pi$. Since $R_{f(p)}\in f(p)=Tr_{\frak{F}}(G)={\cal F}$, we have that $R_{f(p)}\leq R_{\cal F}.$ Hence by the isomorphism $R/R_{f(p)}/R_{\cal F}/R_{f(p)}\simeq R/R_{\cal F},$
we obtain that $R/R_{\cal F}\in \bigcap _{p\in \pi}\frak{E}_{p'}=\frak{E}_{\pi '}.$ Consequently, $R\in {\cal F}\circ \frak{E}_{\pi '}={\cal F}$ by Lemma 4.2.
Therefore ${\cal F}=SLR(f)$.  Note that $G_{\frak{F}}=G_{\cal F}$. Let ${\cal X}=(1),$ then by Lemma 4.4, $G$ is $\frak{F}$-constrained, that is, $C_G(G_{\frak{F}})\subseteq G_{\frak{F}}$.

The following known result directly follows from Lemma 4.4.

{\bf Corollary 4.6.} If $G$ is a soluble group, then $C_G(F(G))\leq F(G).$

{\bf Lemma 4.7.} Let ${\cal F}$ be a $\pi$-semilocal Fitting set of an ${\cal F}\circ\frak{S}^{\pi}$-group $G$ defined by a full and invariable $H$-function $f$. If $V$ is an ${\cal F}$-subgroup of $G$ and $V\geq G_{\cal F}$, then $V/G_f\in \frak{E}_{\pi'}.$

{\bf Proof.} Since $f$ is a full and invariable $H$-function, we may let $f(p)={\cal X}$ for all $p\in \pi.$
By Proposition  3.2, ${\cal X}\subseteq {\cal X}\circ \frak{E}_{p'}$. Hence $f(p)\subseteq \bigcap _{p\in \pi}f(p)\circ \frak{E}_{p'}={\cal F}$. This means that $f$ is also a integrated H-function of ${\cal F}$. Since $G_{\cal F}\unlhd V$ and $G_{\cal X}\leq G_{\cal F},$ $G_{\cal X}=(G_{\cal F})_{\cal X}=G_{\cal F}\cap V_{\cal X}$ by Lemma 2.1. It follows that $G_{\cal X}\leq V_{\cal X}.$ Hence $[V_{\cal X}, G_{\cal F}]\leq V_{\cal X}\cap G_{\cal F}=G_{\cal X}=G_{f(p)}$, and so $V_{\cal X}\leq C_G(G_{\cal F}/G_{\cal X})\leq G_{\cal F}$ by Lemma 4.4. It follows that $V_{\cal X}\leq G_{\cal F}\cap V_{\cal X}=G_{\cal X}$. Consequently, $V_{\cal X}=G_{\cal X}$, that is, $V_{f(p)}=G_{f(p)}$ for all $p\in \pi$. Hence $V_f=G_f.$ It follows from $V\in {\cal F}$ that $V/G_f\in \bigcap _{p\in \pi}\frak{E}_{p'}=\frak{E}_{\pi '}.$

From Lemmas 4.2 and 4.7, we obtain  the following

{\bf Proposition 4.8.} Support that ${\cal F}$ be a $\pi$-semilocal Fitting set of an ${\cal F}\circ \frak{S}^{\pi}$-group $G$ defined by a full and invariable $H$-function $f$, and $V$ a subgroup of $G$ containing $G_{\cal F}.$ Then $V\in {\cal F}$ if and only if $V/G_f\in\frak{E}_{\pi '}.$

\section{Proof and Some Applications of Theorem A}

{\bf Proof of Theorem A.} (1) Support that $G\in {\cal F}\circ \frak{S}^{\pi}$, where $\pi=\sigma ({\cal F})$.

Firstly assume that $2\in \pi$. Then every $\pi '$-factor of $G/G_{\cal F}$ is soluble by the well known Feit-Thompson Theorem. But as $G/G_{\cal F}$ is $\pi$-soluble, it follows that $G/G_{\cal F}$ is soluble. Hence by Lemma 2.5, $G$ possesses an ${\cal F}$-injector and any two ${\cal F}$-injectors are conjugate.

Now assume that $2\notin \pi$. Since $\pi=\sigma ({\cal F})$ and $G_{\cal F}\in {\cal F}$,  $G_{\cal F}$ is a $\pi$-group. Hence  $G_{\cal F}$ is soluble by the Feit-Thompson theorem. But as $G/G_{\cal F}\in \frak{S}^{\pi}$, we have that $G$ is $\pi$-soluble. It follows from Lemma 2.3  that (1) holds.

(2) Support that ${\cal F}$ is $\pi$-semilocal and $G$ is $\pi$-soluble. We prove that $G$ possesses  exactly one conjugacy class of  ${\cal F}$-injectors and the index of every ${\cal F}$-injector in $G$ is a $\pi$-number.

Suppose that it is false and let $G$ be a counterexample of minimal order. Then clear $G\neq 1$.

If $G\in \frak{E}_{\pi '}$, then by Lemma 4.2, $G\in \frak{E}_{\pi '}\subseteq {\cal F}\circ \frak{E}_{\pi '}={\cal F}$. In this case, (2) holds.
We may, therefore, assume that $G\notin \frak{E}_{\pi '}.$

Let $M$ be a maximal normal subgroup of $G$. Then the choice of $G$ implies that $M$ possesses exactly one conjugacy class of  ${\cal F}$-injectors and $|M:V|$ is a $\pi$-number for each ${\cal F}$-injector $V$ of $M$. We now prove (2) via the following steps.

(a) {\sl There exists an ${\cal F}$-maximal subgroup of $G$ containing $V$ and the index of it in $G$ is a $\pi$-number.}

By  Lemma 2.2(2), $V$ is an ${\cal F}$-maximal subgroup of $M$. Hence there exists an ${\cal F}$-maximal subgroup $F$ of $G$ containing $V$. Then $V\leq F\cap M\unlhd F$ and $F\cap M\in {\cal F}.$
It follows that $V=F\cap M$ and $|M : F\cap M|$ is a $\pi$-number. Since $G$ is $\pi$-soluble, $G/M$ is either a soluble $\pi$-group or a $\pi '$-group.

Assume that $|G:M|$ is a $\pi$-number. Then $|G:V|=|G:M||M:V|$ is a $\pi$-number. Consequently, $|G:F|$ is a $\pi$-number.

Now assume that $|G:M|$ is a $\pi '$-number.
Since $V=F\cap M$, by Lemma 2.2(4) $(F\cap M)^x$ is also an $\frak{F}$-injector of $M$ for any $x\in G$.  Hence $F\cap M$ and $(F\cap M)^x$ are conjugate in $M$. By Frattini argument,  $G=MN_G(F\cap M)=MN_G(V).$ Hence $|G|=|M||N_G(V)|/|N_M(V)|$.
Since $V\leq M\cap N_G(V)=N_M(V),$  $|M: N_M(V)| <|M:V|$ is a $\pi$-number. Thus we have that $|G:N_G(V)|$ is a $\pi$-number. Since $G$ is $\pi$-soluble,  $G$ possesses exactly one conjugacy class of  Hall $\pi '$--injectors. Hence there exists a Hall $\pi'$-subgroup $H$ such that $H\leq N_G(V)$. Then $HV=F_1$ is a subgroup of $G$ and $|G:F_1|$ is a $\pi$-number. Clearly, $V\unlhd F_1$. Hence $V\leq (F_1)_{\cal F}.$  Moreover, $|F_1 :V|=|H:H\cap V|$ is $\pi '$-number. Thus $F_1/V\in \frak{E}_{\pi '}$. From the isomorphism $F_1/(F_1)_{\cal F}\simeq (F_1/V)/((F_1)_{\cal F}/V)$, it follows that $F_1/(F_1)_{\cal F}\in\frak{E}_{\pi '}.$ In view of Lemma 4.2, $F_1\in {\cal F}\circ \frak{E}_{\pi '}={\cal F}.$  Then by $M\lhd G$, we have that $F_1\cap M\in {\cal F}$, and so $V=F_1\cap M$ since $V$ is an ${\cal F}$-maximal subgroup of $M$. In order to prove that (a), we only need to prove that $F_1$ is an ${\cal F}$-maximal subgroup of $G$.

Assume that $F_1\leq F_2$ and $F_2$ is an ${\cal F}$-maximal subgroup of $G$. Since $G/M\in \frak{E}_{\pi '},$ $G=HM=F_1M$. Hence $F_2=F_1(F_2\cap M).$  As $V=F_1\cap M$ is an ${\cal F}$-injector of $M$, $F_1\cap M$ is an ${\cal F}$-maximal subgroup of $M$. But $F_1\cap M\leq F_2\cap M$ and $F_2\cap M\in {\cal F},$ so $F_1\cap M=F_2\cap M$. Thus $F_1=F_2.$ This shows that (a) holds.

(b) {\sl $G$ possesses an ${\cal F}$-injector and the index of the ${\cal F}$-injector in $G$ is a $\pi$-number. }

Since $G\notin \frak{E}_{\pi'}$, $1\neq  G^{\frak{E}_{\pi'}}.$ We consider the following two possible cases:

{\sl Case 1.  $G^{\frak{E}_{\pi'}} <G$.}

Then there exists a maximal normal subgroup $N$ of $G$ such that $G^{\frak{E}_{\pi'}}\leq N$. By the choice of $G$, $N$ possesses an ${\cal F}$-injectors $V_1$, any two ${\cal F}$-injectors of $N$ are conjugate in $N$ and $|N:V_1|$ is a $\pi$-number. Since $(G/G^{\frak{E}_{\pi '}})/(N/G^{\frak{E}^{\pi '}})\simeq G/N$, $G/N\in \frak{E}_{\pi '}$. By (a),  there exists an ${\cal F}$-maximal subgroup $E$ of $G$ containing $V_1$ and  $|G:E|$ is a $\pi$-number.
Note that if we can prove that $E\cap K$ is an ${\cal F}$-injector subgroup of $K$ for any maximal normal subgroup $K$ of $G$, then $E$ is an ${\cal F}$-injector of $G$ by Lemma 2.2(3).

Let $K$ be a maximal normal subgroup of $G$. By the choice of $G$, $K$ possesses an ${\cal F}$-injector, $V_2$ say, and $|K:V_2|$ is a $\pi$-number.  Let $T$ be an ${\cal F}$-maximal subgroup of $G$ containing $V_2.$  By the choice of $G$ again, $N\cap K$ possesses an ${\cal F}$-injector and any two ${\cal F}$-injectors are conjugate in $N\cap K$. By Lemma 2.2(5), $E_1=V_1\cap (N\cap K)$ is an ${\cal F}$-injector of $N\cap K$. Hence, $E_1$ is an ${\cal F}$-maximal subgroup of $N\cap K$. Since $V_1\cap (N\cap K)\leq E\cap (N\cap K)$ and  $E\cap (N\cap K)\in {\cal F}$ (since $E\in {\cal F}$ and ${\cal F}$ is a Fitting set), so $E_1=V_1\cap (N\cap K) =E\cap (N\cap K)$. Now since $V_1\leq N$, $E_1=V_1\cap K$. Let $E_2=V_2\cap (N\cap K).$  With a similar argument, we have that $E_2=V_2\cap (N\cap K)=T\cap (N\cap K)=V_2\cap N$ and $E_2$ is an ${\cal F}$-injector of $N\cap K$. Hence $E_1$ and $E_2$ are conjugate in $N\cap K$. Without loss of generality, we may assume that $E_1=E_2=F^*$.

Since $G$ is $\pi$-soluble,  $G/K$ is either a soluble $\pi$-group or a $\pi '$-group.

Assume that $G/K\in\frak{S}_{\pi}$. Since $G/N\in\frak{E}_{\pi '}$, $G=G_{\pi '}N$ for some Hall $\pi'$-subgroup $G_{\pi'}$ of $G$. Moreover, as $|G:E|$ is a $\pi$-number, there exists $x\in G$ such that $(G_{\pi'})^x\leq E$. Hence $G=(G_{\pi'})^xN= EN.$  If $G$ has a unique maximal normal subgroup $K$, then $K=N$, so $E_2=V_2\cap K$ is an ${\cal F}$-injector of $K$. This implies that $V_2$ is ${\cal F}$-injector of $G$ by Lemma 2.2(3) and $|G:V_2|=|G:K||K:V_2|$ is a $\pi$-number. Hence, we may assume that $G$ has two different maximal normal subgroup $N$ and $K$.  Then $G=NK,$ so $G=EN=NK.$ By $NK/N\simeq K/N\cap K$ and $NK/K\simeq N/N\cap K$, we have that $K/N\cap K\in \frak{E}_{\pi'}$ and $N/N\cap K\in \frak{S}_{\pi}.$ Then $K=K_{\pi'}(K\cap N)$ for some Hall $\pi'$-subgroup $K_{\pi'}$ of $K$. But as $|G:E|$ is a $\pi$-number, we have that $|KE:E|=|K:K\cap E|$ is a $\pi$-number. Without loss of generality, we may assume that $K_{\pi'}\leq K\cap E$. Then
$$K=(K\cap E)(K\cap N). \eqno(*)$$

We now prove that $V_2$ and $E\cap K$ are conjugate in $K$.

Note that  $V_2$ is an ${\cal F}$-injector of $K$ and $F^*=V_2\cap N=V_2\cap (N\cap K)$ is ${\cal F}$-maximal in $N\cap K$.
By Lemma 2.4, we only need to prove that $E\cap K$ is an ${\cal F}$-maximal in $K$.

Assume that there exists an ${\cal F}$-maximal subgroup $E^*$ of $G$ such that  $E\cap K\leq E^*\leq K.$ Since $F^*=V_1\cap K=V_2\cap N=V_1\cap (K\cap N)=V_2\cap (K\cap N)$ is an ${\cal F}$-injector of $K\cap N$ and $V_1\cap (K\cap N)\leq E\cap (K\cap N)\leq E^*\cap (K\cap N)\in{\cal F}$, we have that $V_1\cap (K\cap N)=E\cap (K\cap N)=E^*\cap (N\cap K).$
Then by the equation (*),
$$E^*=E^*\cap K=E^*\cap (K\cap E)(K\cap N)=(E\cap (N\cap K))(E\cap K)=E\cap K.$$
 This shows that $E\cap K$ is an ${\cal F}$-maximal subgroup of $K$. Thus $V_2$ and $E\cap K$ are conjugate in $K$, so $E\cap K$ is an ${\cal F}$-injector of $K$ for any maximal normal subgroup $K$ of $G$. Thus by Lemma 2.2(3), $E$ is an ${\cal F}$-injector of $G$.

Now assume that $G/K$ is a $\pi'$-group, that is, $G/K\in \frak{E}_{\pi'}$. Since $G/N$ is also a $\pi'$-group, $G/(K\cap N)\in \frak{E}_{\pi'}$. If $K\cap N=1$, then $G\in \frak{E}_{\pi'}$, it is impossible. Hence $K\cap N\neq 1.$ Since $F^*=T\cap (N\cap K)=E\cap (N\cap K)$ is an ${\cal F}$-maximal subgroup of $N\cap K$ and $F^*\leq E\cap T$, $E$ and $T$ are conjugate in $G$ by Lemma 2.4, that is, $E=T^x$ for some $x\in G$.  But by the choice of $G$, $K$ possesses exactly one conjugacy class of ${\cal F}$-injectors. Since $T\cap K\lhd T\in{\cal F}$, we have that $T\cap K\in {\cal F}$. Moreover as $V_2\leq T\cap K\leq K$ and $V_2$ is an ${\cal F}$-maximal in $K$, $V_2=T\cap K$ is an ${\cal F}$-injector of $K$.  Hence $E\cap K =T^x\cap K$ is an ${\cal F}$-injector of $K$ for any maximal normal subgroup $K$ of $G$. It follows from Lemma 2.2(3) that $E$ is an ${\cal F}$-injector of $G$. Thus (b) holds in the case 1.

{\sl Case 2. $G^{\frak{E}_{\pi'}}=G$.}

In this case, for any maximal normal subgroup $M^*$ of $G$, we have that $G/M^*$ is a $\pi$-group (In fact, if $G/M^*$ is a $\pi'$-group, then $G^{\frak{E}_{\pi'}}\leq M^*$, a contradiction),  so $G/M^*$ is abelian $p$-soluble for some prime $p\in \pi$. Let $M_1$ be a maximal normal subgroup of $G$. The choice of $G$ implies that $M_1$ possesses an ${\cal F}$-injector, say $V_1$. Support that $F_1$ is an ${\cal F}$-maximal subgroup of $G$ containing $V_1$. If we can prove that $F_1\cap M_2$ is an ${\cal F}$-injector of $M_2$ for any maximal normal subgroup $M_2$ of $G$, then $F_1$ is an ${\cal F}$-injector of $G$ by Lemma 2.2(3). In fact, by the choice of $G$,  $M_2$ has an ${\cal F}$-injector $V_2$. Let $F_2$ be an ${\cal F}$-maximal subgroup of $G$ containing $V_2.$ If $M_1\cap M_2=1$, then $G=G/(M_1\cap M_2)\in \frak{S}_{\pi}$. Hence by  Fischer, Gasch$\ddot{u}$tz and Harley Theorem (see \cite[VIII, (2.9)]{DH}), $G$ has an ${\cal F}$-injector, and the index of every ${\cal F}$-injector in $G$ is a $\pi$-number.

Now assume that $M_1\cap M_2\neq 1$. Then $V_1\cap (M_1\cap M_2)$ and $V_2\cap (M_1\cap M_2)$ are ${\cal F}$-injector of  $M_1\cap M_2$ by Lemma 2.2(5), and $V_1\cap (M_1\cap M_2)$ and $V_2\cap (M_1\cap M_2)$ are conjugate in  $M_1\cap M_2$  by the choice of $G$. Without loss of generality, we may assume that $V_1\cap (M_1\cap M_2)= V_2\cap (M_1\cap M_2)=V.$ Then $V\leq F_1\cap F_2$ and $V$ is an ${\cal F}$-maximal subgroup of $M_1\cap M_2$. Note that $G/(M_1\cap M_2)\in \frak{N}_{\pi}.$ Hence by Lemma 2.4, $F_1$ and $F_2$ are conjugate in $G$, that is, $F_1^x=F_2$ for some $x\in G$. As $F_2\cap M_2\lhd F_2\in {\cal F}$, we have that $F_2\cap M_2\in {\cal F}$. But since $V_2$ is ${\cal F}$-maximal in $M_2$, $(F_1\cap M_2)^x=F_1^x\cap M_2=F_2\cap M_2=V_2$. Hence $F_1^x\cap M_2=F_2\cap M_2$ is an ${\cal F}$-injector of $M_2.$  Now by Lemma 2.2(3), $F_2$ is an ${\cal F}$-injector of $G$.

Final, by the choice of $G$, we have that $|M_2 : V_2|$ is $\pi$-number. Then $|G: V_2|=|G: M_2||M_2 :V_2|$ is a $\pi$-number. Thus, (b) holds in the case 2.

{\sl (c) Any two ${\cal F}$-injectors of $G$ are conjugate in $G$.}

Support that $V_1$ and $V_2$ are any two ${\cal F}$-injectors of $G$, and $M$ is a maximal normal subgroup of $G$. Since $G$ is $\pi$-soluble, $G/M$ is either a nilpotent $\pi$-group or a $\pi'$-group. By Lemma 2.2(5), $V_1\cap M$ and $V_2\cap M$ are ${\cal F}$-injectors of $M$, and by the choice of $G$ they are conjugate in $M$. So $V_1^m\cap M=V_2\cap M$ for some $m\in M$. Let $V=V_1^m\cap M=V_2\cap M$. Since $V_1^m$ and $V_2$ are ${\cal F}$-maximal subgroups of $G$ and $V$ is an ${\cal F}$-maximal subgroup of $M$. Hence by Lemma 2.4,$V_1^m$ and $V_2$ are conjugate in $G$. Hence (c) holds.

This completes the proof of (2).

(3) Support that ${\cal F}$ is semilocal and $G/G_{\cal F}$ is $\pi$-soluble, where $\sigma(G_{\cal F})\subseteq \pi.$

Assume that $2\notin \pi$. Then $G_{\cal F}$ is soluble, so $G$ is $\pi$-soluble. In this case, the statement (3) of the theorem follows from (2). Now assume that $2\in \pi$. Then $G/G_{\cal F}$ is soluble by the Feit-Thompson theorem. Hence by Lemma 2.5, $G$ possesses an ${\cal F}$-injector and any two ${\cal F}$-injectors are conjugate.

The theorem is proved.

{\bf Corollary 5.1.} Let $G$ be a $\pi$-soluble group and $\frak{X}$ a nonempty Fitting class such that $\frak{X}\frak{E}_{\pi'}=\frak{X}.$ Then  $G$ possesses exactly one conjugacy class of $\frak{X}$-injectors.

{\bf Proof.} Let ${\cal X}=Tr_{\frak{X}}(G)=Tr_{\frak{X}\frak{E}_{\pi'}}(G)$ and $H\leq G$.  Since $H_{\frak{X}}=H_{\cal X}$, $Tr_{\frak{X}\frak{E}_{\pi'}}(G)=\{H\leq G : H/H_{\frak{X}}\in\frak{E}_{\pi'}\}=\{H\leq G : H/H_{\cal X}\in \frak{E}_{\pi'}\}.$ This implies that ${\cal X}={\cal X}\circ \frak{E}_{\pi'}$,  so ${\cal X}$ is a semilocal by Lemma 4.2. Hence by Theorem A(2), $G$ possesses exactly one conjugacy class of ${\cal X}$-injectors. But since the $\frak{X}$-injectors of $G$ and ${\cal X}$-injector coincide, Thus the corollary holds.

Note that the class of $\pi$-closed groups $\frak{E}_{\pi}\frak{E}_{\pi '}$ is a Fitting class such that $\frak{E}_{\pi}\frak{E}_{\pi '}=(\frak{E}_{\pi}\frak{E}_{\pi '})\frak{E}_{\pi '}$, and the class of $\pi$-special groups $\frak{N}_{\pi}\frak{E}_{\pi '}$ is also a Fitting class such that $\frak{N}_{\pi}\frak{E}_{\pi '}=(\frak{N}_{\pi}\frak{E}_{\pi '})\frak{E}_{\pi '}.$  Hence the following two results directly follows from Corollary 5.1.

{\bf Corollary 5.2.}  Every $\pi$-soluble group $G$ possesses exactly one conjugacy class of $\pi$-closed injectors.

{\bf Corollary 5.3.} Every $\pi$-soluble group $G$ possesses exactly one conjugacy class of $\pi$-special injectors.

Let $\pi =\sigma(G)$, then the following result is  directly follows from Theorem A(2) and Lemma 4.2 in addition to some results in \cite{DH, Fis1, BalE, Shem1, Semen, Bolado, Guo2} mentioned in section 1 can be obtained directly from theorem A.

{\bf Corollary 5.4} \cite[Theorem 3]{And}. Let ${\cal F}$ be a Fitting set of a soluble group $G$. Then $G$ possesses exactly one conjugacy class of ${\cal F}$-injectors.

{\bf Remark 5.5.} In connection with Theorem A(2), if ${\cal F}$ is $\pi$-semilocal but $G\notin {\frak S}^{\pi}$ and suppose that $G$ possesses an ${\cal F}$-injector, then the index of the ${\cal F}$-injector of $G$  in $G$ may be a  $\pi$-number. In fact, we have the following

{\bf Proposition 5.6.} Support that ${\cal F}$ be a $\pi$-semilocal Fitting set of $G$ and $G\in ({\cal F}\circ \frak{S}_{\pi})\cap E_{\pi'},$ where $E_{\pi'}$ is the class of the groups which has a Hall ${\pi'}$-subgroup. Let $V$ be a ${\cal F}$-injector of $G$. Then:

(1) If $K\unlhd G$, then $G=KN_G(V\cap K);$

(2) $|G:V|$ is a $\pi$-number.

{\bf Proof.} (1) Since $K\unlhd G$, $V\cap K$ is a ${\cal F}$-injector of $K$ and $(V\cap K)^x$ is an ${\cal F}$-injector of $K$ for any $x\in G$ by Lemma 2.2(4)(5). By Proposition 3.1, ${\cal F}\circ \frak{S}_{\pi}$ is a Fitting set. As $K\unlhd G\in {\cal F}\circ \frak{S}_{\pi}$, we have that  $K\in {\cal F}\circ \frak{S}$. Hence any two ${\cal F}$-injectors of $K$ are conjugate in $K$ by Theorem A(1), that is, there exists $k\in K$ such that $xk\in N_G(V\cap K).$ It follows that $x\in KN_G(V\cap K)$. Hence $G=KN_G(V\cap K).$

(2) Assume that it is false and let $G$ be a counterexample of minimal order. Then clearly $G\notin {\cal F}$, that is,  $G_{\cal F} < G$.  Let $M$ be a maximal normal subgroup of $G$ containing $G_{\cal F}$. By Lemma 2.2(5), $V\cap M$ is an ${\cal F}$-injector of $M$. The choice of $G$ implies that $|VM:V|=|M:V\cap M|$ is a $\pi$-number. If $G=VM$, then the statement (2) holds, a contradiction. We may, therefore, assume that $G>VM.$ Since $G_{\cal F}\leq M$ and $G/G_{\cal F}$ is soluble, $|G:M|=p$, for some prime $p\in \pi$. Hence $VM=M$, and so $V\leq M$. This implies that  $|M:V|$ is a $\pi$-number. Consequently, $|M:N_G(V)\cap M|=|M:N_M(V)|$ is a $\pi$-number. By (1), $G=MN_G(M\cap V)=MN_G(V).$  Thus $|G|=|M||N_G(V)|/|M\cap N_G(V)|$, and so $|G:N_G(V)|=|M:N_M(V)|$ is a $\pi$-number. But as $G\in E_{\pi'}$, $G$ has a Hall $\pi'$-subgroup. It follows from  $G/G_{\cal F}\in \frak{S}_{\pi}$ that $G_{\pi '}\leq G_{\cal F}\leq V\leq N_G(V).$ Hence $V\unlhd VG_{\pi '}$. Consequently, $V\leq (VG_{\pi'})_{\cal F}$. Since $VG_{\pi'}/V$ is a $\pi '$-group, by the isomorphism
$$(VG_{\pi'}/V)/((VG_{\pi'})_{\cal F}/V)\simeq VG_{\pi'}/(VG_{\pi'})_{\cal F},$$
we obtain that  $VG_{\pi'}/(VG_{\pi'})_{\cal F}\in \frak{E}_{\pi'}.$ It follows from Lemma 4.2 that $VG_{\pi'}\in {\cal F}\circ\frak{E}_{\pi'}={\cal F}.$ But since $V$ is a maximal ${\cal F}$-subgroup of $G$, $VG_{\pi'}=V$ and so $G_{\pi'}\leq V.$  Therefore $|G:V|$ is a $\pi$-number.

\section{Proof and Some Applications of Theorem B}

{\bf Proof of Theorem B.}  Support that ${\cal F}$ be a $\pi$-semilocal Fitting set of $G$ defined by a full H-function $f$  such that $f(p)={\cal X}$ for all $p\in \pi$, for some Fitting set ${\cal X}$ of $G$. Then $G_f=G_{\cal X}.$  Since $G/G_f$ is $\pi$-soluble, $G\in{\cal X}\circ\frak{S}^{\pi}$.
By Proposition 3.2, ${\cal X}\subseteq {\cal X}\circ \frak{E}_{p'}$ for all $p\in \pi$. Hence ${\cal X}\subseteq \bigcap_{p\in \pi}{\cal X}\circ \frak{E}_{p'}={\cal F}$, and so $f$ is integrated.
By the isomorphism $(G/G_{\cal X})/(G_{\cal F}/G_{\cal X})\simeq G/G_{\cal F}$, we have that $G/G_{\cal F}$ is $\pi$-soluble, that is, $G\in{\cal F}\circ \frak{S}^{\pi}$.

(1) Let $V$ is an ${\cal F}$-injector of $G$. Then by Lemma 2.2(1)(2), $V$ is a ${\cal F}$-maximal subgroup of $G$ and $G_{\cal F}\leq V$. By Lemma 4.7, $V/G_{\cal X}\in \frak{E}_{\pi'}.$
Note that a $\pi$-soluble group possesses a Hall $\pi'$-subgroup and any two Hall $\pi'$-subgroups are conjugate (see \cite[Theorem 1.7.7]{Guo1}). Assume that $V/G_{\cal X}<F/G_{\cal X}$ for some Hall $\pi'$-subgroup $F/G_{\cal X}$ of $G/G_{\cal X}$. Then by Lemma 4.3, $F\in {\cal F}$. But $V$ is a maximal ${\cal F}$-subgroup of $G$, so $F=V$. Thus $V/G_f=V/G_{\cal X}$ is a Hall $\pi'$-subgroup of $G/G_{\cal X}.$

Conversely, assume that $V/G_{\cal X}$ is a Hall $\pi'$-subgroup of $G/G_{\cal X}$. Let $M$ be any maximal normal subgroup of $G$. Then $M_{\cal X}=G_{\cal X}\cap M$ by Lemma 2.1(1).

If $G_{\cal X}\leq M$, then $M_{\cal X}=G_{\cal X}.$  Hence $(V\cap M)/M_{\cal X}$ is a Hall $\pi'$-subgroup of $M/M_{\cal X}$. By induction, $V\cap M$ is an ${\cal F}$-injector of $M$. As $V/G_{\cal X}\in\frak{E}_{\pi'}$, $V\in {\cal F}$ by Lemma 4.3. Assume that $V\leq W$ for some ${\cal F}$-maximal subgroup of $G$. Note that $G_{\cal F}/G_{\cal X}=G_{\cal F}/(G_{\cal F})_{\cal X}=G_{\cal F}/(G_{\cal F})_{f(p)}\in \bigcap _{p\in \pi}\frak{E}_{p'}=\frak{E}_{\pi'}$.  Hence $G_{\cal F}/G_{\cal X}\leq V/G_{\cal X}$, so $G_{\cal F}\leq W$. Then by Lemma 4.7, $W/G_{\cal X}\in \frak{E}_{\pi'}.$  It follows that $V/G_{\cal X}=W/G_{\cal X}$ and so $V=W$. Therefore by Lemma 2.2(3), we obtain that $V$ is an ${\cal F}$-injector of $G$.

Now assume that $G_{\cal X}\nleq M$. Then $G=G_{\cal X}M$ and so $G/G_{\cal X}\simeq M/M\cap G_{\cal X}=M/M_{\cal X}$. Since $V/G_{\cal X}$ is a Hall $\pi'$-subgroup of $G/G_{\cal X}$,
$(V\cap M)/M_{\cal X}$ is a Hall $\pi'$-subgroup of $M/M_{\cal X}$. By induction, $V\cap M$ is an ${\cal F}$-injector of $M$. With a similar argument as above, we have that $V$ is an ${\cal F}$-injector of $G$.

Thus the statement (1) of Theorem B holds.

(2) Since $G/G_f$ is $\pi$-soluble, $G/G_f$ possesses a Hall $\pi'$-subgroup and any two Hall $\pi'$-subgroups are conjugate in $G$. Hence (2) follows from (1).

(3) If $V$ is an ${\cal F}$-maximal subgroup of $G$ and $G_{\cal F}\leq V$. Then by Lemma 4.7, $V/G_f\in\frak{E}_{\pi'}$. Assume that $V/G_f < F/G_f$ for some Hall $\pi'$-subgroup $F/G_f$ of $G/G_f$. Since $G_{\cal F}\leq V <F$,  $F\in {\cal F}$ by Proposition 4.8. Hence $F=V$. Consequently, $V/G_f$ is a Hall $\pi'$-subgroup of $G/G_f$. Then by (1), $V$ is an ${\cal F}$-injecotr of $G$. The converse directly follows from the definition of ${\cal F}$-injector of $G$.

(4) Let $V$ be an ${\cal F}$-injector of $G$. Then by (1), $V/G_f$ is a Hall $\pi'$-subgroup of $G/G_f$. Since $G$ is $\pi$-soluble, there exists a Hall $\pi'$-subgroup $G_{\pi'}$ of $G$ such that $V/G_f=G_{\pi'}G_f/G_f$. Hence $V=G_{\pi'}G_f$, and so (4) holds.

{\bf Corollary 6.1.} The $\pi$-closed injectors of a $\pi$-soluble group $G$ are just the subgroups of the type $G_{\pi'}O_{\pi}(G)$. Hence the $\pi$-closed injectors are the maximal $\pi$-closed subgroups of $G$ containing the $\pi'$-closed radical $G_{\pi '}$ of $G$.

{\bf Proof.} Let $\frak{F}$ be the class of all $\pi$-closed groups, that is, $\frak{F}=\frak{E}_{\pi}\frak{E}_{\pi'}$.  Since $G$ is $\pi$-soluble, every Hall $\pi$-subgroup of $G$ is soluble. Let ${\cal F}=Tr_{\frak{E}_{\pi}\frak{E}_{\pi'}}(G)$. Then ${\cal F}=Tr_{\frak{S}_{\pi}\frak{E}_{\pi'}}(G)$. If we can prove that ${\cal F}$ is a $\pi$-semilocal Fitting set of $G$ defined by a full and invariable H-function $f$, then by Theorem B(4), every ${\cal F}$-injector of $G$ is of the type $G_{\pi '}G_f$, where $G_{\pi '}$ is some Hall $\pi '$-subgroup of $G$.

Let ${\cal X}=Tr_{\frak{S}_{\pi}}(G)$ and $f(p)={\cal X}$ for all $p\in \pi$. Then $f$ is invariable.
We now prove that $f$ is a full H-function of ${\cal F}$.
Clearly, ${\cal X}\subseteq {\cal X}\circ \frak{N}_p$ for any $p\in \pi$. Let $H\leq G$ and $H\in {\cal X}\circ \frak{N}_p$. Then $H/H_{\cal X}\in \frak{N}_p\subseteq \frak{S}_{\pi}$ and $H_{\cal X}\in{\cal X}\subseteq \frak{S}_{\pi}.$ Hence $H\in \frak{S}_{\pi}$, so $H\in {\cal X}.$ Thus ${\cal X}\circ \frak{N}_p={\cal X}$, that is, $f$ is a full H-function.

We now prove that ${\cal F}=SLR(f)=\bigcap_{p\in \pi}{\cal X}\circ \frak{E}_{p'}={\cal X}\circ \frak{E}_{\pi'}$.

If $H\in {\cal F}$, then $H\in \frak{S}_{\pi}\frak{E}_{\pi'}$ and so $H/H_{\frak{S}_{\pi}}\in\frak{E}_{\pi'}.$ Since $H_{\frak{S}_{\pi}}=H_{\cal X}$, $H\in {\cal X}\circ\frak{E}_{\pi'}$.
This shows that ${\cal F}\subseteq {\cal X}\circ\frak{E}_{\pi'}.$  Conversely, assume that $L\leq G$ and $L\in {\cal X}\circ\frak{E}_{\pi'}.$ Then $L/L_{\cal X}\in\frak{E}_{\pi'}.$
As $L_{\cal X}=L_{\frak{S}_{\pi}},$ we have that $L\in \frak{S}_{\pi}\frak{E}_{\pi'}.$ Hence $L\in {\cal F}$. This shows that ${\cal F}={\cal X}\circ\frak{E}_{\pi'}.$  This completes the proof.

{\bf Corollary 6.2.} Let $G$ be a $\pi$-soluble group and $\frak{F}_k$ the class of all $\pi$-soluble groups with $\pi$-length $\leq k$, where $k$ is nonnegative integer. Then the $\frak{F}_k$-injectors of $G$ are the $\frak{F}_k$-maximal subgroups containing $G_{{\cal F}_k}.$

{\bf Proof.} Let ${\cal F}_k=Tr_{\frak{F}_k}(G)$. If $k=0$, then ${\cal F}_k=Tr_{\frak{E}_{\pi'}}(G)$ and the ${\cal F}_k$-injectors are Hall $\pi'$-subgroups of $G$. Hence the statement of the corollary hold. Now assume that $k>0.$  With a similar argument as in the proof of Corollary 6.1, we can prove that ${\cal F}_k$ is a semilocal Fitting set of $G$ defined by a full and invariable H-function $f$ such that $f(p)=Tr_{\frak{F}_{k-1}\frak{S}_{\pi}}(G)$ for all $p\in \pi$, where $\frak{F}_{k-1}$ is the Fitting class of all $\pi$-soluble groups with $\pi$-length $\leq k-1$. Hence the statement of the corollary follows from Theorem B(3).

\end{document}